
\documentstyle{amsppt}
\TagsOnRight

\hoffset=.1in

\define\ducha{\chi_{{\phantom{}}_\a}} 
\define\duchbe{\chi_{{\phantom{}}_\be}} 
\define\duchde{\chi_{{\phantom{}}_\de}}

\define\adx{\operatorname{ad}}

\define\a{\alpha}
\define\be{\beta}
\define\de{\delta}
\define\vf{\varphi}
\define\lm{\lambda}
\define\ve{\varepsilon}

\define\si{\sigma}
\define\om{\omega}

\define\vd{\varDelta}
\define\vg{\varGamma}

\define\caa{\Cal A}
\define\car{\Cal R}

\define\ek{\sum_k}

\define\0{Hei\-sen\-berg group}
\define\1{Hei\-sen\-berg algebra}
\define\2{quantum Lie algebra}
\define\3{coproduct}
\define\4{equation}
\define\5{functional}
\define\6{calculation}
\define\7{calculate}
\define\8{structure}
\define\9{infinitesimal}

\define\dix{dimensional}

\define\wor{Wo\-ro\-no\-wicz}

\define\fol{following}
\define\bim{bimodule}
\define\linv{left-invariant}
\define\rinv{right-invariant}
\define\cou{counterpart}

\define\xx{\,\otimes\,}
\define\pri{\pi r^{-1}}

\topmatter
\title Differential calculus\\
on the quantum Hei\-sen\-berg group
\endtitle
\rightheadtext{Hei\-sen\-berg group}
\author Piotr Kosi\'nski\\
{\it Department of Theoretical Physics}\\
{\it  University of \L \'od\'z}\\
{\it ul. Pomorska 149/153, 90--236 \L \'od\'z, Poland}\\
Pawe\l \/ Ma\'slanka$^*$, Karol Przanowski $^*$\\
{\it Department of Functional Analysis}\\
{\it  University of \L \'od\'z}\\
{\it ul. St. Banacha 22, 90--238, \L \'od\'z, Poland}
\endauthor
\leftheadtext{P. Kosi\'nski, P. Ma\'slanka, K. Przanowski}
\thanks
*\ \ \ Supported by \L \'od\'z University grant N$^o$ 505$\slash$445
\endthanks

\abstract The differential calculus on the quantum Hei\-sen\-berg group is 
con-\linebreak structed. The duality between quantum Hei\-sen\-berg group 
and algebra is proved.
\endabstract

\endtopmatter
\document

\head I. Introduction
\endhead

The one \dix{} deformed \0 and algebra were investigated in [1], [2]. In 
this paper, using \wor{}'s theory ([3]), we construct the differential 
calculus on the deformed one \dix{} \0 and we describe the \8 of its \2. 
Then we prove that our \2 is equivalent to the one \dix{} deformed \1.

\head II. The differential calculus 
\endhead

The quantum group $H(1)_q$ is a matrix quantum group \`a la \wor{} ([4])
$$
T = \pmatrix
1 & \a & \be\\
0 & 1 & \de\\
0 & 0 & 1
\endpmatrix
\tag{1}
$$
where the matrix elements $\a$, $\be$, $\de$ generate the algebra $\caa$ and 
satisfy the \fol{} relations ([1])
$$
\aligned
& [\a,\be] = i\lm \a,\\
& [\de,\be] = i\lm \de,\\
& [\a,\de] = 0,
\endaligned
\tag{2}
$$
$\lm$ being a real parameter.

The \3, counit and antipode are given by
$$
\aligned
& \vd(\a) =  I \xx \a + \a \xx I,\\
& \vd(\be) =  I \xx \be + \be \xx I + \a \xx \de,\\
& \vd(\de) =  I \xx \de + \de \xx I,\\
& S(\a) = - \a,\\
& S(\be) =  - \be + \a  \de,\\
& S(\de) =  - \de,\\
& \ve(\a) = \ve(\be) = \ve(\de) = 0.
\endaligned
\tag{3}
$$
The main ingredient of the \wor{} theory is the choice of a right ideal in 
$\ker \ve$, which is invariant under the adjoint action of the group. The 
adjoint action is defined as follows
$$
\adx(a) = \ek b_k \xx S(a_k) c_k
\tag{4}
$$
here
$$
(\vd \xx I) \circ \vd(a) = \ek a_k \xx b_k \xx c_k.
$$

One can prove the following

\proclaim{Theorem 1}  Let $\car \subset \ker \ve$ be the right ideal 
generated by the  \fol{} elements: $\a^2$, $\de^2$, $\be\a$, $\be\de$, 
$\a\de$, $\be^2 + 2i\lm \be$. Then
\roster
\item"{(i)}" $\car$ is $\adx$-invariant,  $\adx(\car) \subset \car \xx \caa$
\item"{(ii)}" $ \ker \ve\slash\,\car$ is spanned by the \fol{} elements: 
$\a$, $\be$, $\de$.
\endroster
\endproclaim 

Having established the \8 of $\car$ we follow closely the \wor{} 
construction. The basis of the space of the \linv{} 1-forms consists of the 
\fol{} elements 
$$
\aligned
& \om_\a \equiv \pri (I \xx \a) = d\a ,\\
& \om_\be \equiv \pri (I \xx \be) = d\be,\\
& \om_\de \equiv \pri (I \xx \de) =  d\be - \a d\de ;
\endaligned
\tag{5}
$$
here the mapping $r^{-1}$ is given by
$$
r^{-1}  (a \xx b) = (a \xx I) (S \xx I)\vd (b), \qquad \qquad a,b \in \caa,
$$
and the mapping $\pi$ is given by
$$
\pi(\ek a_k \xx b_k) = \ek a_k d b_k
$$
where $\ek a_k \xx b_k \in \caa \xx \caa$ is such an element that
$$
\ek a_k b_k = 0.
$$
The next step is to find the commutation rules between the invariant forms 
and generators of $\caa$. The detailed \6s result in the \fol{} formulae
$$
\aligned
& [\a,\om_\a] = 0,\\
& [\de,\om_\a] = 0,\\
& [\be,\om_\a] = -i\lm \om_\a,\\
& [\a,\om_\de] = 0,\\
& [\de,\om_\de] = 0,\\
& [\be,\om_\de] = -i\lm \om_\de,\\
& [\a,\om_\be] = 0,\\
& [\de,\om_\be] = 0,\\
& [\be,\om_\be] = 2i\lm \om_\be.
\endaligned
\tag{6}
$$
Then, following \wor{}'s paper [3], we can  construct the \rinv{} forms
$$
\aligned
& \eta_\a = \om_\a,\\
& \eta_\de = \om_\de,\\
& \eta_\be = \om_\be - \om_\a \de + \om_\de \a.
\endaligned
\tag{7}
$$
This concludes the description of the \bim{} $\vg$ of $1$-forms on $H(1)_q$. 
The external algebra can now be constructed as follows ([3]). On $\vg^{\xx 
2}$ we define a \bim{} homomorphism $\si$ such that
$$
\si (\om \xx_{{}_{\!\!\!\caa}} \ \eta) = \eta \xx_{{}_{\!\!\!\caa}} \ \om
\tag{8}
$$
for any \linv{} $\om \in \vg$ and  any \rinv{} $\eta \in \vg$. Then by 
definition
$$
\vg^{\wedge 2}  = \dfrac{\vg^{\xx 2}}{\ker(I - \si)}.
\tag{9}
$$
Equations (7)--(9) allow us to \7 the external product of \linv{} $1$-forms. The result reads
$$
\aligned
& \om_\be  \wedge  \om_\a = - \om_\a \wedge \om_\be,\\
& \om_\be  \wedge  \om_\de = - \om_\de \wedge \om_\be,\\
& \om_\be  \wedge  \om_\be = 0,\\
& \om_\a  \wedge \om_\a = 0,\\
& \om_\de  \wedge \om_\de = 0,\\
& \om_\a  \wedge \om_\de = -\om_\de \wedge  \om_\a.
\endaligned
\tag{10} 
$$
To complete the external calculus, we derive the Cartan-Maurer \4s
$$
\aligned
&  d\om_\a = 0,\\
& d \om_\de = 0,\\
& d\om_\be = - \om_\de \wedge \om_\be.
\endaligned
\tag{11} 
$$
\vskip.3cm

\head III. Quantum Lie algebra
\endhead

In order to obtain the \cou{} of the classical Lie algebra, we introduce the 
\cou{} of the \linv{} vector fields. They are defined by the formula
$$
da = (\ducha \ast a)\om_\a + (\duchbe \ast a)\om_\be + (\duchde \ast 
a)\om_\de.
\tag{12}
$$
In order to find the \2, we apply the external derivative to both sides of 
(12), we use $d^2a = 0$ on the left-hand side and \7 the right-hand side 
using (11) and again (12). Nullifying the coefficients in front of basis 
elements of $\vg^{\wedge 2}$, we find the \2
$$
\aligned
& [\ducha, \duchbe] = 0,\\
& [\duchde, \duchbe] = 0,\\
& [\ducha  \duchde] = \duchbe.
\endaligned
\tag{13}
$$
From the \wor{} theory, it follows that the coproduct of the \5 $\vf_i$ 
($\vf_i \equiv  \ducha, \duchbe, \duchde)$ can be written in the form
$$
\vd \vf_i = \sum_j \vf_j \xx f_{ji} + I \xx \vf_i
\tag{14}
$$
where $f_{ji}$ are the \5s entering in the commutation rules between the 
\linv{} forms and elements of $\caa$
$$
\om_j a = \sum_i (f_{ji} \ast a)\om_i.
\tag{15}
$$
Then, it follows from commutation rules (6) that the coproduct for our \5s 
can be written in the \fol{} form
$$
\aligned
& \vd \ducha = \ducha  \xx f_\a + I \xx \ducha ,\\
& \vd \duchbe =  \duchbe \xx f_\be + I \xx \duchbe ,\\
& \vd \duchde = \duchde \xx f_\de + I \xx \duchde
\endaligned
\tag{16}
$$
Using the fact that $\vd f_i = f_i \xx f_i$ ($i = \a,\be,\de$) ([3]) and (6) 
and (15), we can \7 the \5s $f_i$. After some \6s we obtain
$$
\aligned
& f_\a = (I - 2i\lm \duchbe)^{\frac{1}{2}},\\
& f_\be = I - 2i\lm \duchbe,\\
& f_\de = (I - 2i\lm \duchbe)^{\frac{1}{2}}.
\endaligned
\tag{17}
$$
Now it is easy to see that the substitution
$$
\aligned
& \duchde = B_0,\\
& \duchbe = B_1,\\
& \ducha = B_2
\endaligned
\tag{18}
$$
reproduces the \8 of the Hopf algebra generated by the \9 generators 
(obtained by contraction procedure) of the quantum matrix pseudogroup 
$H(1)_q$, which was described in [1]. This proves the duality between the 
quantum \0 and algebra.

\enddocument